\documentclass[10pt]{article}
\usepackage{amsmath,amssymb,amsthm,latexsym}
\usepackage{graphicx,xcolor,tikz,enumerate}
\renewenvironment{proof}{\noindent{\bf Proof.}}{~~$\Box$}
\theoremstyle{plain}
\theoremstyle{thmstyleone}%
\newtheorem{thrm}{Theorem}
\newtheorem{prot}[thrm]{Proposition}%

\theoremstyle{thmstyletwo}%
\newtheorem{exa}{Example}%
\newtheorem{remk}{Remark}%

\theoremstyle{thmstylethree}%
\newtheorem{dfnt}{Definition}%
\begin{document}
\title{On the existence of internal uninorms on a bounded lattice: the antichain case
}
\author{
Hua-Peng Zhang\quad Yao Ouyang\thanks{Corresponding author.\quad Email:huapengzhang@163.com(H.-P. Zhang); oyy@huznu.edu.cn(Y. Ouyang)}\\
\small\it {School of Science, Nanjing University of Posts and Telecommunications, Nanjing, Jiangsu 210023, China}\\
\small\it Faculty of Science, Huzhou Normal University, Huzhou, Zhejiang 313000, China}

\date{}
\maketitle
\begin{abstract}
\hspace{4mm}
In the context of bounded lattices, internal uninorms are ordered semigroups with the distinctive property that their binary operation always yields one of its two arguments. Here we focus on the case where the set of elements that are incomparable with a prescribed neutral element is an antichain¡ªthat is, its members are all pairwise incomparable. Within this framework, we first characterize the precise structure of internal uninorms possessing such a neutral element. Based on this structural result, we then establish a necessary and sufficient criterion for their existence, formulated purely in the language of the lattice order.

\end{abstract}

\section{Introduction}\label{sec1}

Algebraically, triangular norms (t-norms) and triangular conorms (t-conorms)~\cite{SS60} on the unit interval are canonical examples of ordered commutative monoids, with neutral elements $1$ and $0$, respectively. Beyond their mathematical elegance, t-norms have become indispensable in applied domains such as fuzzy logic~\cite{EG01} and information fusion~\cite{HGGG11}. The pivotal generalization to uninorms, proposed by Yager and Rybalov~\cite{YR96}, relaxed the fixed-neutral-element constraint by permitting a neutral element $e$ anywhere within $[0,1]$. This unifying framework has given rise to a rich taxonomy of subclasses, among which $\mathcal{U}{\min}$, $\mathcal{U}{\max}$, and representable uninorms have been studied extensively. The theory of uninorms is well-summarized in the survey \cite{SZMMD25}.

The study of uninorms on bounded lattices was initiated by Kara\c{c}al and Mesiar~\cite{KM15}, who extended this class of aggregation functions from the real unit interval to the general context of bounded lattices. A fundamental result of their work is the construction of the smallest and greatest uninorms for any prescribed neutral element. Following this foundational extension, research activity on lattice-based uninorms has flourished along two primary trajectories. The first trajectory concerns the development of diverse construction techniques, leading to a growing repository of uninorm families (e.g.,~\cite{J21, OZ20, XZ23}). Concurrently, a distinct and vital trajectory seeks to provide axiomatic characterizations for significant subclasses. Exemplary works in this vein include the characterization of the $\mathcal{U}{\min}$ and $\mathcal{U}{\max}$ classes by Zhang et al.~\cite{ZWWOD20}, and the characterization of idempotent uninorms on specific lattice types by Mesiarov\'{a}-Zem\'{a}nkov\'{a} et al.~\cite{MMSZ251, MMSZ252}.


An internal uninorm~\cite{CKM19} is an idempotent uninorm with the additional property that its output is always one of its two input values. It is known that every bounded lattice admits an idempotent uninorm with any chosen neutral element. However, the existence of an internal uninorm is not guaranteed, making the characterization of such existence a natural and important problem. In recent research, Ouyang et al.~\cite{OZWD24} solved this problem under the assumption that the set $I_e$ of elements incomparable with the neutral element $e$ is a chain. They derived a complete structural representation and established the following precise existence condition: either all elements in $I_e$ have a common meet with $e$, or they have a common join with $e$. Equivalently, for any two elements in $I_e$, the one with strictly smaller meet with $e$ must have a join with $e$ that is strictly larger than or equal to the join of the other with $e$.

However, when all elements incomparable with the neutral element $e$ do not constitute a chain, the existence of internal uninorms has not been characterized yet. The purpose of this paper is to consider the existence problem in the case that all elements incomparable with $e$ constitute an antichain. We first characterize the structure of internal uninorms in this case. Based on the structural characterization theorem for internal uninorms, we explore the existence of internal uninorms. As we shall see, the condition that for any two elements in $I_e$, the one with strictly smaller meet with $e$ must have a join with $e$ that is strictly larger than or equal to the join of the other with $e$, is also necessary and sufficient for the existence of internal uninorms in the antichain case.

\section{Preliminaries}\label{sec2}

This section outlines the basic order-theoretic notions and terminology employed throughout the paper. We refer to the classic text by Birkhoff~\cite{B73} for further background on lattice theory.

A \emph{partially ordered set} (poset) $(P,\leq)$ consists of a nonempty set $P$ equipped with a reflexive, antisymmetric and transitive binary relation $\leq$. If $P$ contains a greatest element (denoted $1$) and a least element (denoted $0$), we call it a \emph{bounded poset}. For $a,b\in P$, we write $a<b$ when $a\leq b$ and $a\neq b$. If neither $a\leq b$ nor $b\leq a$ holds, we say $a$ and $b$ are \emph{incomparable} and write $a\parallel b$. The set of all elements incomparable with a fixed $a\in P$ is denoted by $I_{a}$. A poset in which every pair of elements is comparable is called a \emph{chain}; if every pair of distinct elements in a subset is incomparable, that subset is called an \emph{antichain}.

For $a,b\in P$ with $a\leq b$, we use the standard notation $[a,b]$, $]a,b]$, $[a,b[$ and $]a,b[$ for closed, half-open, and open intervals. A subset $A\subseteq P$ is an \emph{upper set} (respectively, \emph{lower set}) if $x\in A$ and $x\leq y$ (respectively, $y\leq x$) imply $y\in A$~\cite{GHKLMS03}.

A \emph{lattice} is a poset $(L,\leq)$ in which any two elements $x,y$ have a greatest lower bound, called the \emph{meet} $x\wedge y$, and a least upper bound, called the \emph{join} $x\vee y$. Every chain is a lattice; in that case, meet and join coincide with the minimum and maximum with respect to the total order.

A binary operation $F:X\times X\to X$ on a nonempty set $X$ is said to be \emph{internal} (or conservative)~\cite{CouDevMar19} provided $F(x,y)\in\{x,y\}$ for all $x,y\in X$. Internal operations are necessarily idempotent. The following structural result characterises associative internal operations.

\begin{thrm}[\cite{CouDevMar19}]\label{thrm-neworder}
A binary operation $F$ on a nonempty set $X$ is commutative, associative and internal if and only if there exists a total order $\preccurlyeq$ on $X$ such that $F=\max_\preccurlyeq$.
\end{thrm}

\begin{remk}
By the Order-Extension Principle, any partial order can be extended to a total order; consequently, on any nonempty set there exists at least one commutative, associative and internal binary operation.
\end{remk}

In what follows, $L$ always stands for a bounded lattice $(L,\leq,\wedge,\vee,0,1)$. Aggregation functions on such lattices-in particular, uninorms-are defined as follows.

\begin{dfnt}[\cite{KM15}]\label{def:uninorm}
A binary operation $U:L^2\to L$ is a \emph{uninorm} with neutral element $e\in L$ if it satisfies, for all $x,y,z,w\in L$:
\begin{enumerate}[{\rm(i)}]
    \item \emph{Increasingness}: if $x\leq y$ and $z\leq w$, then $U(x,z)\leq U(y,w)$;
    \item \emph{Commutativity}: $U(x,y)=U(y,x)$;
    \item \emph{Associativity}: $U(U(x,y),z)=U(x,U(y,z))$;
    \item \emph{Neutral element}: $U(e,x)=U(x,e)=x$.
\end{enumerate}
\end{dfnt}

\noindent
Commutativity makes the increasingness property equivalent to monotonicity in each variable separately. Together, increasingness and associativity endow $L$ with the structure of an ordered semigroup.

An internal uninorm is always idempotent. The converse, however, does not hold in general; idempotent uninorms need not be internal. On a \emph{chain}, the two notions coincide: every idempotent uninorm is internal (see~\cite[Proposition 2.3]{OZWD22}).

We denote by $\mathcal{U}_{\mathrm{Int}}(L,e)$ the set of all internal uninorms on $L$ with neutral element $e$. When this set is nonempty, the subset $L\setminus I_e$ necessarily forms a chain under the induced order~\cite{CKM19}. Hence, for any $U\in\mathcal{U}_{\mathrm{Int}}(L,e)$, its restriction to $(L\setminus I_e)^2$ is an idempotent uninorm on the bounded chain $L\setminus I_e$. Since idempotent uninorms on bounded chains are already fully described~\cite{OZWD22,SMM23}, the remaining task for characterising internal uninorms on a general bounded lattice reduces to:
\begin{itemize}
    \item describing their behaviour on the incomparable part $I_e^2$, and
    \item understanding how this behaviour interacts with the structure on $L\setminus I_e$.
\end{itemize}

A necessary condition for the existence of an internal uninorm was given in~\cite[Proposition 3.3]{OZWD24}:
\begin{equation}\label{eq:nec-cond}
\forall u,v\in I_e:\qquad u\wedge e < v\wedge e \;\Longrightarrow\; u\vee e \geq v\vee e. \tag{$\star$}
\end{equation}
When $I_e$ is a chain, condition~\eqref{eq:nec-cond} is also sufficient~\cite[Theorem 4.11]{OZWD24}. \textbf{The main contribution of the present paper} is to prove that~\eqref{eq:nec-cond} remains sufficient in the opposite extremal case, namely when $I_e$ is an antichain (Theorem~\ref{thrm-chainAndantichain} below). This completes the picture for the two most elementary order-theoretic configurations of the incomparable set $I_e$.

\section{Structure of internal uninorms}\label{sec3}
In this section, we characterize the structure of any internal uninorm on $L$ with neutral element $e$ when $I_e$ is an antichain.

Let $U\in\mathcal U_{\,\rm Int}(L, e)$ and $u\in I_e$. Define $A_u\subseteq[0, e[\,$ and $B_u\subseteq\,]e, 1]$ as follows:
\begin{equation}\label{eq-Au}
A_u=\{x\in[0, e[\ \mid U(x, u)=x\}\,;
\end{equation}
\begin{equation}\label{eq-Bu}
B_u=\{x\in\,]e, 1]\mid U(x, u)=x\}\,.
\end{equation}
The following proposition collects some facts about $\{A_u\}_{u\in I_e}$ and $\{B_u\}_{u\in I_e}$ as well as some properties of $U$ related to $\{A_u\}_{u\in I_e}$ and $\{B_u\}_{u\in I_e}$.

\begin{prot}\cite{OZWD24}\label{Prot-ProAuBu}
Let $U\in\mathcal U_{\,\rm Int}(L, e)$ and $u, v\in I_e$. Then the following results hold:
\begin{enumerate}[\rm(i)]
\item $A_u\subseteq[0, u\wedge e]$ is a lower set and $B_u\subseteq[u\vee e, 1]$ is an upper set;

\item $v\wedge e\in A_u$ or $u\vee e\in B_v$;

\item $A_u=[0, u\wedge e]$ or $B_u=[u\vee e, 1]$;

\item $U(x, y)=y$ for all $(x, y)\in([0, e[\,\setminus A_u)\times [u, 1]\cup (\,]e, 1]\setminus B_u)\times [0, u]$;

\item $U(u, v)=v$ implies that $A_v\subseteq A_u$ and $B_v\subseteq B_u$;

\item $u\parallel v$ and $U(u, v)=v$ imply that $v\wedge e\in  A_u$ and $v\vee e\in B_u$;

\item $u\parallel v$ and $U(u, v)=v$ imply that $v\wedge e\leq u\wedge e$ and $v\vee e\geq u\vee e$.
\end{enumerate}
\end{prot}
\begin{proof}
(i)--(v) come from Proposition 3.2 of~\cite{OZWD24}.

(vi)  If $v\wedge e\notin A_u$, then the increasingness of $U$ implies that
\[
u=U(u, v\wedge e)\leq U(u, v)=v,
\]
a contradiction with $u\parallel v$. Similarly, we have $v\vee e \in B_u$.

(vii) It follows from (vi) and (i).
\end{proof}

When $I_e$ is an antichain, the structure of any internal uninorm on $L$ with neutral element $e$ can be characterized as follows.
\begin{thrm}\label{thrm-characterizationthrm}
Let $L\setminus I_e$ be a chain, $I_e$ be an antichain and $U$ be a binary operation on $L$. Then $U\in {\mathcal U}_{\rm Int}(L, e)$ if and only if there exist a family of lower sets $\{{\mathcal A}_u\}_{u\in I_e}$ with ${\mathcal A}_u\subseteq [0, u\wedge e]$, a family of upper sets $\{{\mathcal B}_u\}_{u\in I_e}$ with ${\mathcal B}_u\subseteq [u\vee e, 1]$, an internal uninorm $U_0$ on $L\setminus I_e$ with neutral element $e$ and a commutative, associative and internal binary operation $Q$ on $I_e$ such that
\begin{eqnarray}\label{eq-UniInt}
  U(x, y) =\left \{
        \begin {array}{ll}
         U_0(x, y)
                  &\quad {\rm if} \ \ (x, y)\in (L\setminus I_e)^2  \\[1mm]
         x
                  &\quad {\rm if} \ \ (x, y)\in I_e\times (([0, e]\setminus\mathcal A_x)\cup (\,]e, 1]\setminus\mathcal B_x))\bigcup (\mathcal A_y\cup\mathcal B_y)\times I_e \\[1mm]
         y
                  &\quad {\rm if} \ \ (x, y)\in I_e\times (\mathcal A_x\cup\mathcal B_x)\bigcup (([0, e]\setminus\mathcal A_y)\cup (\,]e, 1]\setminus\mathcal B_y))\times I_e \\[1mm]
        Q(x, y)
                  &\quad {\rm if} \ \ (x, y)\in I_e^2
        \end {array}
       \right.
\end{eqnarray}
and the following properties hold for all $u, v\in I_e$:
\begin{itemize}
\item[\rm(C1)] $U_0(x, y)=y$ for all $(x, y)\in([0, e[\,\setminus\mathcal A_u)\times [u\vee e, 1]\cup(\,]e, 1]\setminus\mathcal B_u)\times [0, u\wedge e]$;

\item[\rm(C2)] $u\not=v$ and $Q(u, v)=v$ imply that $v\wedge e\in {\mathcal A}_u$ and $v\vee e\in {\mathcal B}_u$.
\end{itemize}
\end{thrm}
\begin{proof}
{\it Necessity:} Suppose that $U\in\mathcal U_{\,\rm Int}(L, e)$. For every $u\in I_e$, define $\mathcal A_u=A_u$ and $\mathcal B_u=B_u$, where $A_u$ and $B_u$ are given by~\eqref{eq-Au} and~\eqref{eq-Bu}, respectively. Let $U_0=U|_{(L\setminus I_e)^2}$ and $Q=U|_{I_e^2}$ be the restrictions of $U$ to $(L\setminus I_e)^2$ and $I_e^2$, respectively. Clearly, $U_0$ is an internal uninorm on $L\setminus I_e$ with neutral element $e$, $Q$ is a commutative, associative and internal binary operation on $I_e$, and $U$ is representable as in~\eqref{eq-UniInt}. It follows from Proposition~\ref{Prot-ProAuBu}(i) that ${\mathcal A}_u\subseteq [0, u\wedge e]$ is a lower set and ${\mathcal B}_u\subseteq [u\vee e, 1]$ is an upper set. (C1) and (C2) follow from Proposition \ref{Prot-ProAuBu}(iv) and (vi), respectively.

{\it Sufficiency:}  Obviously, we only need to prove that $U$ is increasing and associative.

Increasingness: Let $x, y, z\in L$ with $y<z$, we can verify that $U(x, y)\leq U(x, z)$. Since $U$ is internal, the inequality trivially holds if $x=y$ or $x=z$. We consider other cases. Note that $\{y, z\}\nsubseteq I_e$ as $I_e$ is an antichain.
\begin{enumerate}[\rm(i)]
\item $y\in I_e$ and $z\in [y\vee e, 1]$.

If $x\in I_e$, $U(x, y)=x$, and $z\in\mathcal B_x$, then $U(x, y)<x\vee e\leq z=U(x, z)$.

If $x\in I_e$, $U(x, y)=x$, and $z\notin\mathcal B_x$, then $U(x, y)=x=U(x, z)$.

If $x\in I_e$ and $U(x, y)=y$, then it follows from (C2) that $y\vee e\in {\mathcal B}_x$, which implies that $z\in {\mathcal B}_x$, and thus $U(x, y)<y\vee e\leq z=U(x, z)$.

If $x\in [e, 1]$, then $U(x, y)\leq x\vee y\leq x\vee z=U(x, z)$.

If $x\in [0, e[\,\cap\mathcal A_y$, then $U(x, y)=x=x\wedge z\leq U(x, z)$.

If $x\in [0, e[\,\setminus\mathcal A_y$, then it follows from (C1) that $U(x, z)=z$, and thus $U(x, y)=y<z=U(x, z)$.

\item $z\in I_e$ and $y\in [0, z\wedge e]$.

The proof is analogous to that of case (i).

\item $y, z\in L\setminus I_e$.

If $x\in L\setminus I_e$, then $U(x, y)=U_0(x, y)\leq U_0(x, z)=U(x, z)$.

If $x\in I_e$ and $y\in\mathcal A_x$, then $U(x, y)=y\leq x\wedge z\leq U(x, z)$.

If $x\in I_e$ and $z\in \mathcal B_x$, then $U(x, y)\leq x\vee y\leq z=U(x, z)$.

If $x\in I_e$, $y\notin\mathcal A_x$ and $z\notin \mathcal B_x$, then $z\notin\mathcal A_x$ and $y\notin \mathcal B_x$, and thus
$U(x, y)=x=U(x, z)$.
\end{enumerate}

Associativity: Let $x, y, z\in L$. We can verify that $U(U(x, y), z)=U(x, U(y, z))$. Since $U_0$ and $Q$ are associative, the equality naturally holds if $\{x, y, z\}\subseteq L\setminus I_e$ or $\{x, y, z\}\subseteq I_e$.
We consider the remaining cases.
\begin{enumerate}[(i)]
\item Only one of $x, y, z$ belongs to $I_e$.

Since $U$ is commutative, we only need to consider the case that $x\in I_e$ and $y, z\in L\setminus I_e$.

If $y\in \mathcal A_x\cup\mathcal B_x$, then it is easy to prove that $U_0(y, z)\in \mathcal A_x\cup\mathcal B_x$ by the aid of (C1), and thus $U(U(x, y), z)=U_0(y, z)=U(x, U(y, z))$.

If $y\notin \mathcal A_x\cup\mathcal B_x$ and $z\notin \mathcal A_x\cup\mathcal B_x$, then $U(y, z)\notin \mathcal A_x\cup\mathcal B_x$, which implies that $U(U(x, y), z)=x=U(x, U(y, z))$.

If $y\notin \mathcal A_x\cup\mathcal B_x$ and $z\in \mathcal A_x\cup\mathcal B_x$, then it is easy to prove that $U_0(y, z)=z$ by the aid of (C1), and thus $U(U(x, y), z)=U(x, z)=U(x, U(y, z))$.

\item Exactly two of $x, y, z$ belong to $I_e$.

Since $U$ is commutative, we only need to consider the case that $x, y\in I_e$ and $z\in L\setminus I_e$.

If $z\in \mathcal A_x\cup\mathcal B_x$ and $z\in \mathcal A_y\cup\mathcal B_y$, then $z\in \mathcal A_{U(x, y)}\cup\mathcal B_{U(x, y)}$, and thus
$U(U(x, y), z)=z=U(x, z)=U(x, U(y, z))$.

If $z\in \mathcal A_x\cup\mathcal B_x$ and $z\notin \mathcal A_y\cup\mathcal B_y$, then it is easy to prove that $U(x, y)=y$ by (C2), and thus $U(U(x, y), z)=U(x, y)=U(x, U(y, z))$.

If $z\notin \mathcal A_x\cup\mathcal B_x$ and $z\in \mathcal A_y\cup\mathcal B_y$, then $U(x, y)=x$ by (C2), and thus
$U(U(x, y), z)=U(x, z)=U(x, U(y, z))$.

If $z\notin \mathcal A_x\cup\mathcal B_x$ and $z\notin \mathcal A_y\cup\mathcal B_y$, then $z\notin \mathcal A_{U(x, y)}\cup\mathcal B_{U(x, y)}$, and thus $U(U(x, y), z)=U(x, y)=U(x, U(y, z))$.
\end{enumerate}
We conclude that $U\in {\mathcal U}_{\rm Int}(L, e)$.
\end{proof}

The following two examples show that (C1) and (C2) in Theorem~\ref{thrm-characterizationthrm} are mutually independent.
\begin{exa}\label{exa-independenceofC2}
Consider the bounded lattice $L=\{0, a, b, e, f, 1, u, v\}$ shown in Figure~\ref{figure-2}.
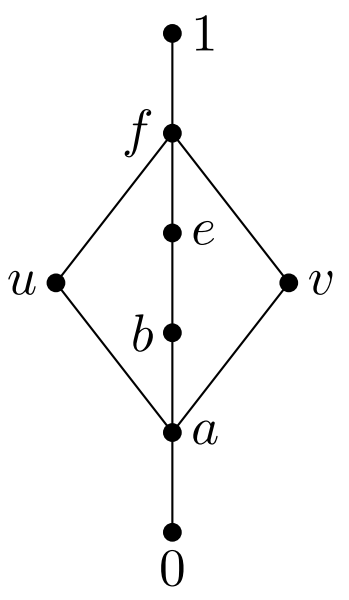
\begin{figure}
\centering
\begin{tikzpicture}[scale = 1.1]
  \draw ( 0, 10) -- (0, 13);
  \node [below] at (0, 10) {$0$};
  \draw[fill] ( 0, 10) circle [radius = 0.05];
  \node [right] at (0, 10.6) {$a$};
  \draw[fill] ( 0, 10.6) circle [radius = 0.05];
  \node [left] at (0, 11.2) {$b$};
  \draw[fill] ( 0, 11.2) circle [radius = 0.05];
  \node [right] at (0, 11.8) {$e$};
  \draw[fill] ( 0, 11.8) circle [radius = 0.05];
  \node [left] at (0, 12.4) {$f$};
  \draw[fill] ( 0, 12.4) circle [radius = 0.05];
  \node [right] at (0, 13.0) {$1$};
  \draw[fill] ( 0, 13.0) circle [radius = 0.05];
  \draw (0, 10.6)--(0.7, 11.5);
  \draw (0, 12.4)--(0.7, 11.5);
  \draw (0, 10.6)--(-0.7, 11.5);
  \draw (0, 12.4)--(-0.7, 11.5);
  \node [left] at (-0.7, 11.5) {$u$};
  \draw[fill] (-0.7, 11.5) circle [radius = 0.05];
  \node [right] at (0.7, 11.5) {$v$};
  \draw[fill] (0.7, 11.5) circle [radius = 0.05];
 \end{tikzpicture}
\caption{Hasse diagram of the lattice in Examples~\ref{exa-independenceofC2} and~\ref{exa-independenceofC1}.}
\label{figure-2}
\end{figure}
Let $U_0$ be the internal uninorm on $L\setminus I_e$ with neutral element $e$ shown in Table~\ref{tab-C2},
\begin{table*}
\centering
\begin{tabular}{c|c c c c c c}
  \hline
  $U_0$ & $0$ & $a$ & $b$ & $e$ & $f$ & $1$ \\ \hline
  $0$ & $0$ & $0$ & $0$ & $0$ & $f$ & $1$ \\
  $a$ & $0$ & $a$ & $a$ & $a$ & $f$ & $1$ \\
  $b$ & $0$ & $a$ & $b$ & $b$ & $f$ & $1$ \\
  $e$ & $0$ & $a$ & $b$ & $e$ & $f$ & $1$ \\
  $f$ & $f$ & $f$ & $f$ & $f$ & $f$ & $1$ \\
  $1$ & $1$ & $1$ & $1$ & $1$ & $1$ & $1$ \\
\end{tabular}
\caption{The internal uninorm $U_0$ in Example~\ref{exa-independenceofC2}.}
\label{tab-C2}
\end{table*}
$Q$ be the commutative, associative and internal binary operation on $I_e$ defined by $Q(u, u)=u$ and $Q(u, v)=Q(v, u)=Q(v, v)=v$, $\mathcal A_u=\mathcal A_v=\{0\}$ and $\mathcal B_u=\mathcal B_v=\{f, 1\}$.
It is easy to check that (C1) holds as $\,]e, 1]\setminus \mathcal B_u=\,]e, 1]\setminus\mathcal B_v=\emptyset$, but (C2) fails as $v\wedge e\notin {\mathcal A}_u$.
\end{exa}

\begin{exa}\label{exa-independenceofC1}
Consider the bounded lattice $L=\{0, a, b, e, f, 1, u, v\}$ shown in Figure~\ref{figure-2}.
Let $U_0$ be the internal uninorm on $L\setminus I_e$ with neutral element $e$ shown in Table~\ref{tab-C1},
\begin{table*}
\centering
\begin{tabular}{c|c c c c c c c c}
  \hline
  $U_0$ & $0$ & $a$ & $b$ & $e$ & $f$ & $1$ \\ \hline
  $0$ & $0$ & $0$ & $0$ & $0$ & $0$ & $0$ \\
  $a$ & $0$ & $a$ & $a$ & $a$ & $a$ & $a$ \\
  $b$ & $0$ & $a$ & $b$ & $b$ & $b$ & $b$ \\
  $e$ & $0$ & $a$ & $b$ & $e$ & $f$ & $1$ \\
  $f$ & $0$ & $a$ & $b$ & $f$ & $f$ & $1$ \\
  $1$ & $0$ & $a$ & $b$ & $1$ & $1$ & $1$ \\
\end{tabular}
\caption{The internal uninorm $U_0$ in Example~\ref{exa-independenceofC1}.}
\label{tab-C1}
\end{table*}
$Q$ be the commutative, associative and internal binary operation on $I_e$ defined by $Q(u, u)=u$ and $Q(u, v)=Q(v, u)=Q(v, v)=v$, $\mathcal A_u=\mathcal A_v=\{0, a\}$ and $\mathcal B_u=\mathcal B_v=\{f, 1\}$.
It is easy to see that (C2) holds, but (C1) fails as $U_0(b, f)=b\not=f$ for $(b, f)\in \bigl([0, e[\,\setminus\mathcal A_u\bigr)\times [u\vee e, 1]$.
\end{exa}

\section{Existence of internal uninorms}\label{sec4}

In this section, we study the existence of an internal uninorm on $L$ with neutral element $e$ when $I_e$ is an antichain.
In view of Theorem~\ref{thrm-characterizationthrm}, the existence of an internal uninorm on $L$ with neutral element $e$ amounts to that of $\{\mathcal A_u\}_{u\in I_e}$, $\{\mathcal B_u\}_{u\in I_e}$, $U_0$ and $Q$ satisfying (C1) and (C2).
First, we have the following observations about $\{\mathcal A_u\}_{u\in I_e}$ and $\{\mathcal B_u\}_{u\in I_e}$.
\begin{prot}\label{prot-ObserAuBu}
Let $\{\mathcal A_u\}_{u\in I_e}$ be a family of lower sets with $\mathcal A_u\subseteq [0, u\wedge e]$, $\{\mathcal B_u\}_{u\in I_e}$ be a family of upper sets with $\mathcal B_u\subseteq [u\vee e, 1]$ and $Q$ be a commutative, associative and internal binary operation on $I_e$ satisfying {\rm(C2)}. Then the following results hold:
\begin{enumerate}[\rm(i)]
\item For any $u\in I_e$, we have $\bigcup\limits_{v\in S(u)}[0, v\wedge e]\subseteq\mathcal A_u$ and $\bigcup\limits_{v\in G(u)}[v\vee e, 1]\subseteq\mathcal B_u$, where $S(u)=\{v\in I_e\mid v\wedge e<u\wedge e\}$ and $G(u)=\{v\in I_e\mid v\vee e>u\vee e\}$;

\item For any $u\in I_e$, if there exists $v\in I_e\setminus\{u\}$ with $u\wedge e=v\wedge e$ and $Q(u, v)=v$, then $\mathcal A_u=[0, u\wedge e]$;

\item For any $u\in I_e$, if there exists $v\in I_e\setminus\{u\}$ with $u\vee e=v\vee e$ and $Q(u, v)=v$, then $\mathcal B_u=[u\vee e, 1]$;

\item For any $u, v\in I_e$, if $u\wedge e=v\wedge e$ and $u\vee e<v\vee e$, then $\mathcal A_u=[0, u\wedge e]$;

\item For any $u, v\in I_e$, if $u\vee e=v\vee e$ and $v\wedge e<u\wedge e$, then $\mathcal B_u=[u\vee e, 1]$;

\item If $u\wedge e=v\wedge e$ for all $u, v\in I_e$, then $\mathcal A_u=[0, u\wedge e]$ for all $u\in I_e\setminus\{u_0\}$, where $u_0$ is the annihilator of $Q$, \emph{i.e.}, $Q(u_0, u)=u_0$ for all $u\in I_e$. In particular, $\mathcal A_u=[0, u\wedge e]$ for all $u\in I_e$ if $Q$ has no annihilator;

\item If $u\vee e=v\vee e$ for all $u, v\in I_e$, then $\mathcal B_u=[u\vee e, 1]$ for all $u\in I_e\setminus\{u_0\}$, where $u_0$ is the annihilator of $Q$. In particular, $\mathcal B_u=[u\vee e, 1]$ for all $u\in I_e$ if $Q$ has no annihilator;

\item If $u\wedge e=v\wedge e$ for all $u, v\in I_e$ and $\{u\vee e\mid u\in I_e\}$ has no greatest element, then $\mathcal A_u=[0, u\wedge e]$ for all $u\in I_e$;

\item If $u\vee e=v\vee e$ for all $u, v\in I_e$ and $\{u\wedge e\mid u\in I_e\}$ has no smallest element, then $\mathcal B_u=[u\vee e, 1]$ for all $u\in I_e$.
\end{enumerate}
\end{prot}
\begin{proof}
(i)--(vii) follow directly from (C2), while (viii) and (ix) follow from (iv) and (v), respectively.
\end{proof}

Whether $I_e$ is an antichain or not, the existence of an internal uninorm on $L\setminus I_e$ with neutral element $e$ fulfilling (C1) is determined by $\{\mathcal A_u\}_{u\in I_e}$ and $\{\mathcal B_u\}_{u\in I_e}$.
\begin{prot}\label{prot-existenceU0}
Let $L\setminus I_e$ be a chain, $\{\mathcal A_u\}_{u\in I_e}$ be a family of lower sets with $\mathcal A_u\subseteq [0, u\wedge e]$ and $\{\mathcal B_u\}_{u\in I_e}$ be a family of upper sets with $\mathcal B_u\subseteq [u\vee e, 1]$. Then the following two statements are equivalent:
\begin{enumerate}[\rm(i)]
\item There exists an internal uninorm $U_0$ on $L\setminus I_e$ with neutral element $e$ fulfilling {\rm (C1)};

\item $v\wedge e\in\mathcal A_u$ or $u\vee e\in\mathcal B_v$ for all $u, v\in I_e$.
\end{enumerate}
\end{prot}
\begin{proof}
The proof is essentially the same as the proof of Proposition 4.8 of~\cite{OZWD24}. For the sake of self-containedness, we provide a succinct proof.

(i) $\Longrightarrow$ (ii) Suppose that there exist $u, v\in I_e$ such that $v\wedge e\notin\mathcal A_u$ and $u\vee e\notin\mathcal B_v$.
Then $(v\wedge e, u\vee e)\in([0, e[\,\setminus\mathcal A_u)\times [u\vee e, 1]$ and
$(u\vee e, v\wedge e)\in(\,]e, 1]\setminus\mathcal B_v)\times [0, v\wedge e]$. Hence, it follows from (C1) that $u\vee e=U_0(v\wedge e, u\vee e)=v\wedge e$, a contradiction.

(ii) $\Longrightarrow$ (i) Define a binary operation $U_0$ on $L\setminus I_e$ as follows:
\begin{eqnarray}\label{eq-ConU0}
U_0(x, y)=\left \{
        \begin {array}{ll}
         x\wedge y,
                  & {\rm if} \ \ (x, y)\in [0, e]^2\cup\bigcup\limits_{u\in I_e}([0, u\wedge e]\times(\,]e, 1]\setminus\mathcal B_u)\cup(\,]e, 1]\setminus\mathcal B_u)\times[0, u\wedge e]) \\[1mm]
         x\vee y,
                  & {\rm otherwise.} \\[1mm]
        \end {array}
       \right.
\end{eqnarray}
Then $U_0$ is an internal uninorm on $L\setminus I_e$ with neutral element $e$. It remains to show that $U_0$ satisfies~(C1), \emph{i.e.}, $U_0(x, y)=y$ for all
\[
(x, y)\in\bigcup\limits_{u\in I_e}(([0, e[\,\setminus\mathcal A_u)\times [u\vee e, 1]\cup(\,]e, 1]\setminus\mathcal B_u)\times [0, u\wedge e])\,.
\]
It is clear that $U_0(x, y)=y$ for all $(x, y)\in\bigcup\limits_{u\in I_e}(\,]e, 1]\setminus\mathcal B_u)\times [0, u\wedge e]$. It follows from (ii) that
\[
(([0, e[\,\setminus\mathcal A_u)\times [u\vee e, 1])\cap([0, v\wedge e]\times(\,]e, 1]\setminus\mathcal B_v))=\emptyset\,,
\]
for all $u, v\in I_e$. Hence, $U_0(x, y)=y$ for all $(x, y)\in\bigcup\limits_{u\in I_e}([0, e[\,\setminus\mathcal A_u)\times [u\vee e, 1]$.
\end{proof}

Whether $I_e$ is an antichain or not, the existence of $Q$ fulfilling (C2) is also determined by $\{\mathcal A_u\}_{u\in I_e}$ and $\{\mathcal B_u\}_{u\in I_e}$.
\begin{prot}\label{prot-existenceQ}
Let $\{\mathcal A_u\}_{u\in I_e}$ be a family of lower sets with $\mathcal A_u\subseteq [0, u\wedge e]$ and $\{\mathcal B_u\}_{u\in I_e}$ be a family of upper sets with $\mathcal B_u\subseteq [u\vee e, 1]$. Then the following two statements are equivalent:
\begin{enumerate}[\rm(i)]
\item There exists a commutative, associative and internal binary operation $Q$ on $I_e$ fulfilling {\rm (C2)};

\item For all $u, v\in I_e$ with $u\not=v$, at least one of the following two conditions holds: {\rm (a)} $v\wedge e\in {\mathcal A}_u$ and $v\vee e\in {\mathcal B}_u$; {\rm (b)} $u\wedge e\in {\mathcal A}_v$ and $u\vee e\in {\mathcal B}_v$.
\end{enumerate}
\end{prot}
\begin{proof}
(i) $\Longrightarrow$ (ii) It follows from the assumption that $Q$ is internal and commutative.

(ii) $\Longrightarrow$ (i)
Define a binary relation $\sim$ on $I_e$ as follows:
\begin{center}
$u\sim v$ if $u=v$, or $u\not=v$ and both (a) and (b) hold.
\end{center}
It is easy to verify that $\sim$ is an equivalence relation on $I_e$. For $u\in I_e$, the equivalence class containing $u$ is defined as $[u]=\{v\in I_e\mid v\sim u\}$. For every equivalence class $[u]$, let $\preccurlyeq_{[u]}$ be an arbitrary total order on $[u]$. Define a binary relation $\preccurlyeq$ on $I_e$ as follows:
\begin{center}
$u\preccurlyeq v$ if $u\sim v$ and $u\preccurlyeq_{[u]} v$, or $u\not\sim v$ and (a) holds.
\end{center}
It is easy to prove that $\preccurlyeq$ is a total order on $I_e$.
Hence, it follows from Theorem~\ref{thrm-neworder} that $Q=\max_\preccurlyeq$ is a commutative, associative and internal binary operation on $I_e$. Clearly, $Q$ fulfills (C2).
\end{proof}

Next, we characterize the existence of $\{\mathcal A_u\}_{u\in I_e}$ and $\{\mathcal B_u\}_{u\in I_e}$ fulfilling Proposition~\ref{prot-existenceU0}(ii) and Proposition~\ref{prot-existenceQ}(ii).
\begin{prot}\label{prot-existenceAuBu}
The following results hold:
\begin{enumerate}[\rm(i)]
\item There exist a family of lower sets $\{\mathcal A_u\}_{u\in I_e}$ with $\mathcal A_u\subseteq [0, u\wedge e]$ and a family of upper sets $\{\mathcal B_u\}_{u\in I_e}$ with $\mathcal B_u\subseteq [u\vee e, 1]$ fulfilling {\rm Proposition~\ref{prot-existenceU0}(ii)} if and only if $u\wedge e<v\wedge e\Longrightarrow u\vee e\geq v\vee e$ for all $u, v\in I_e$;

\item There exist a family of lower sets $\{\mathcal A_u\}_{u\in I_e}$ with $\mathcal A_u\subseteq [0, u\wedge e]$ and a family of upper sets $\{\mathcal B_u\}_{u\in I_e}$ with $\mathcal B_u\subseteq [u\vee e, 1]$ fulfilling {\rm Proposition~\ref{prot-existenceQ}(ii)} if and only if $u\wedge e<v\wedge e\Longrightarrow u\vee e\geq v\vee e$ for all $u, v\in I_e$;

\item There exist a family of lower sets $\{\mathcal A_u\}_{u\in I_e}$ with $\mathcal A_u\subseteq [0, u\wedge e]$ and a family of upper sets $\{\mathcal B_u\}_{u\in I_e}$ with $\mathcal B_u\subseteq [u\vee e, 1]$ fulfilling {\rm Proposition~\ref{prot-existenceU0}(ii)} and {\rm Proposition~\ref{prot-existenceQ}(ii)} if and only if $u\wedge e<v\wedge e\Longrightarrow u\vee e\geq v\vee e$ for all $u, v\in I_e$.
\end{enumerate}
\end{prot}
\begin{proof}
(i) If $u, v\in I_e$ with $u\wedge e<v\wedge e$, then $v\wedge e\notin\mathcal A_u$. It follows from Proposition~\ref{prot-existenceU0}(ii) that $u\vee e\in\mathcal B_v$, and thus $u\vee e\geq v\vee e$.

Conversely, consider $\mathcal A_u=[0, u\wedge e]$ and $\mathcal B_u=[u\vee e, 1]$ for all $u\in I_e$. Since $u\wedge e<v\wedge e\Longrightarrow u\vee e\geq v\vee e$ for all $u, v\in I_e$, it is easy to verify that Proposition~\ref{prot-existenceU0}(ii) holds.

(ii) If $u, v\in I_e$ with $u\wedge e<v\wedge e$, then $u\not=v$ and $v\wedge e\notin\mathcal A_u$. It follows from Proposition~\ref{prot-existenceQ}(ii) that $u\wedge e\in {\mathcal A}_v$ and $u\vee e\in {\mathcal B}_v$, and thus $u\vee e\geq v\vee e$.

Conversely, consider $\mathcal A_u=[0, u\wedge e]$ and $\mathcal B_u=[u\vee e, 1]$ for all $u\in I_e$. Since $u\wedge e<v\wedge e\Longrightarrow u\vee e\geq v\vee e$ for all $u, v\in I_e$, it is easy to verify that Proposition~\ref{prot-existenceQ}(ii) holds.

(iii) The necessity follows from both (i) and (ii), while the sufficiency follows from the proof of (i) and (ii) by
considering $\mathcal A_u=[0, u\wedge e]$ and $\mathcal B_u=[u\vee e, 1]$ for all $u\in I_e$.
\end{proof}

Now, we can obtain a necessary and sufficient condition for the existence of an internal uninorm with a given neutral element.
\begin{thrm}\label{thrm-chainAndantichain}
Let $L\setminus I_e$ be a chain and $I_e$ be an antichain. Then ${\mathcal U}_{\rm Int}(L, e)\neq\emptyset$ if and only if $u\wedge e<v\wedge e\Longrightarrow u\vee e\geq v\vee e$ for all $u, v\in I_e$.
\end{thrm}
\begin{proof}
It follows from Theorem~\ref{thrm-characterizationthrm} and Propositions~\ref{prot-existenceU0},~\ref{prot-existenceQ} and~\ref{prot-existenceAuBu}(iii).
\end{proof}

Under the necessary and sufficient condition, as Theorem~\ref{thrm-characterizationthrm}, and Propositions~\ref{prot-existenceU0} and~\ref{prot-existenceQ} show, to construct $U\in{\mathcal U}_{\rm Int}(L, e)$, we need to construct $\{\mathcal A_u\}_{u\in I_e}$ and $\{\mathcal B_u\}_{u\in I_e}$ fulfilling Proposition~\ref{prot-existenceU0}(ii) and Proposition~\ref{prot-existenceQ}(ii). Generally, Proposition~\ref{prot-ObserAuBu} and the proof of Proposition~\ref{prot-existenceAuBu}(iii) suggest that there seems to be no more choice except considering $\mathcal A_u=[0, u\wedge e]$ and $\mathcal B_u=[u\vee e, 1]$ for all $u\in I_e$. However, there are other options in several particular cases.
\begin{prot}\label{prot-ConsAuBu}
Suppose that $u\wedge e<v\wedge e\Longrightarrow u\vee e\geq v\vee e$ for all $u, v\in I_e$. Then the following results hold:
\begin{enumerate}[\rm(i)]
\item If $u\wedge e\not=v\wedge e$ for all $u, v\in I_e$ with $u\not=v$, then $\{\mathcal A_u\}_{u\in I_e}$ and $\{\mathcal B_u\}_{u\in I_e}$ fulfill {\rm Proposition~\ref{prot-existenceU0}(ii)} and {\rm Proposition~\ref{prot-existenceQ}(ii)}, where $\mathcal A_u=\bigcup\limits_{v\in S(u)}[0, v\wedge e]$ and $\mathcal B_u=[u\vee e, 1]$ for all $u\in I_e$;

\item If $u\vee e\not=v\vee e$ for all $u, v\in I_e$ with $u\not=v$, then $\{\mathcal A_u\}_{u\in I_e}$ and $\{\mathcal B_u\}_{u\in I_e}$ fulfill {\rm Proposition~\ref{prot-existenceU0}(ii)} and {\rm Proposition~\ref{prot-existenceQ}(ii)}, where $\mathcal A_u=[0, u\wedge e]$ and $\mathcal B_u=\bigcup\limits_{v\in G(u)}[v\vee e, 1]$ for all $u\in I_e$.
\end{enumerate}
\end{prot}
\begin{proof}
The proof is omitted as an easy exercise.
\end{proof}

\begin{exa}\label{Exa-1}
Consider the bounded lattice $L=\{0, a, b, e, c, d, 1, u_1, u_2, u_3\}$ shown in Figure~\ref{figure-1}.
\begin{figure}
\centering
\begin{tikzpicture}[scale = 1.1]
  \draw ( 0, 10) -- (0, 13.6);
  \draw[fill] ( 0, 10) circle [radius = 0.05];
  \node [below] at (0, 10) {$0$};
  \draw[fill] ( 0, 10.6) circle [radius = 0.05];
  \node [right] at (0, 10.6) {$a$};
  \draw[fill] ( 0, 11.2) circle [radius = 0.05];
  \node [left] at (0, 11.2) {$b$};
  \draw[fill] ( 0, 11.8) circle [radius = 0.05];
  \node [right] at (0, 11.8) {$e$};
  \draw[fill] (0, 12.4) circle [radius = 0.05];
   \node [left] at (0, 12.4) {$c$};
  \draw[fill] (0, 13.0) circle [radius = 0.05];
   \node [right] at (0, 13.0) {$d$};
  \draw[fill] (0, 13.6) circle [radius = 0.05];
  \node [above] at (0, 13.6) {$1$};
  \draw(0, 11.2)--(0.6, 11.8);
  \draw(0.6, 11.8)--(0, 12.4);
  \draw[fill](0.6, 11.8) circle [radius = 0.05];
  \node [right] at (0.6, 11.8) {$u_2$};
  \draw(0, 10.6)--(1.2, 11.8);
  \draw(1.2, 11.8)--(0, 13.0);
  \draw[fill](1.2, 11.8) circle [radius = 0.05];
  \node [right] at (1.2, 11.8) {$u_1$};
  \draw(0, 10)--(-1.2, 11.8);
  \draw(-1.2, 11.8)--(0, 13.0);
  \draw[fill](-1.2, 11.8) circle [radius = 0.05];
  \node [left] at (-1.2, 11.8) {$u_3$};
\end{tikzpicture}
\caption{Hasse diagram of the lattice $L$ in Example~\ref{Exa-1}.}
\label{figure-1}
\end{figure}
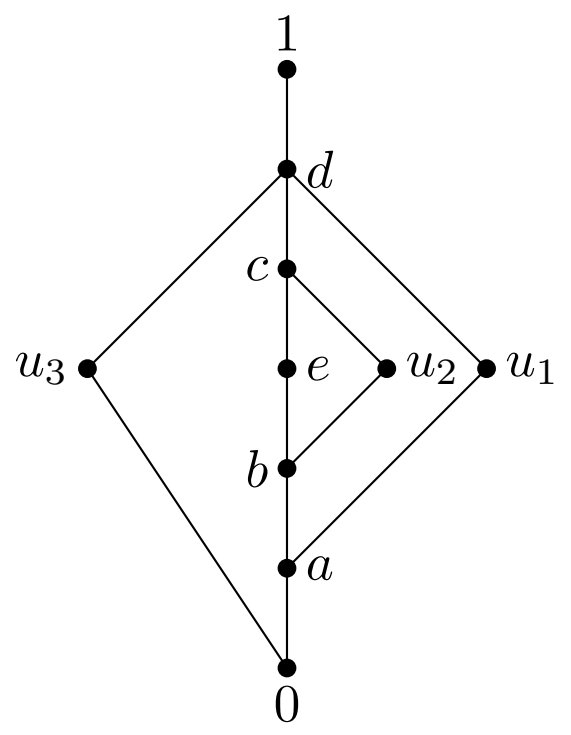
In Example 3.8 of~\cite{OZWD24}, we constructed the internal uninorm $U_1$ on $L$ with neutral element $e$
shown in Table~\ref{tab-IntUni-1}.
\begin{table*}
\centering
\begin{tabular}{c|c c c c c c c c c c}
  \hline
  $U_1$ & $0$ & $a$ & $b$ & $e$ & $c$ & $d$ & $1$ & $u_1$ & $u_2$ & $u_3$ \\ \hline
  $0$ & $0$ & $0$ & $0$ & $0$ & $0$ & $0$ & $0$ & $0$   & $0$   & $0$ \\
  $a$ & $0$ & $a$ & $a$ & $a$ & $a$ & $d$ & $1$ & $u_1$ & $a$   & $u_3$ \\
  $b$ & $0$ & $a$ & $b$ & $b$ & $b$ & $d$ & $1$ & $u_1$ & $b$   & $u_3$ \\
  $e$ & $0$ & $a$ & $b$ & $e$ & $c$ & $d$ & $1$ & $u_1$ & $u_2$ & $u_3$ \\
  $c$ & $0$ & $a$ & $b$ & $c$ & $c$ & $d$ & $1$ & $u_1$ & $c$   & $u_3$ \\
  $d$ & $0$ & $d$ & $d$ & $d$ & $d$ & $d$ & $1$ & $d$   & $d$   & $u_3$ \\
  $1$ & $0$ & $1$ & $1$ & $1$ & $1$ & $1$ & $1$ & $1$   & $1$   & $u_3$ \\
  $u_1$ & $0$ & $u_1$ & $u_1$ & $u_1$ & $u_1$ & $d$ & $1$ & $u_1$  & $u_1$ & $u_3$ \\
  $u_2$ & $0$ & $a$ & $b$ & $u_2$ & $c$ & $d$ & $1$ & $u_1$ & $u_2$ & $u_3$ \\
  $u_3$ & $0$ & $u_3$ & $u_3$ & $u_3$ & $u_3$ & $u_3$ & $u_3$ & $u_3$ & $u_3$ & $u_3$ \\ \hline
\end{tabular}
\caption{The internal uninorm $U_1$ in Example~\ref{Exa-1}.}
\label{tab-IntUni-1}
\end{table*}
Note that $u_1\wedge e\not=u_2\wedge e\not=u_3\wedge e\not=u_1\wedge e$. Due to Proposition~\ref{prot-ConsAuBu}(i), we can construct a different internal uninorm $U_2$ on $L$ with neutral element $e$
shown in Table~\ref{tab-IntUni-2}.
\begin{table*}
\centering
\begin{tabular}{c|c c c c c c c c c c}
  \hline
  $U_2$ & $0$ & $a$ & $b$ & $e$ & $c$ & $d$ & $1$ & $u_1$ & $u_2$ & $u_3$ \\ \hline
  $0$ & $0$ & $0$ & $0$ & $0$ & $0$ & $d$ & $1$ & $0$   & $0$   & $u_3$ \\
  $a$ & $0$ & $a$ & $a$ & $a$ & $a$ & $d$ & $1$ & $u_1$ & $a$   & $u_3$ \\
  $b$ & $0$ & $a$ & $b$ & $b$ & $c$ & $d$ & $1$ & $u_1$ & $u_2$   & $u_3$ \\
  $e$ & $0$ & $a$ & $b$ & $e$ & $c$ & $d$ & $1$ & $u_1$ & $u_2$ & $u_3$ \\
  $c$ & $0$ & $a$ & $c$ & $c$ & $c$ & $d$ & $1$ & $u_1$ & $c$   & $u_3$ \\
  $d$ & $d$ & $d$ & $d$ & $d$ & $d$ & $d$ & $1$ & $d$   & $d$   & $d$ \\
  $1$ & $1$ & $1$ & $1$ & $1$ & $1$ & $1$ & $1$ & $1$   & $1$   & $1$ \\
  $u_1$ & $0$ & $u_1$ & $u_1$ & $u_1$ & $u_1$ & $d$ & $1$ & $u_1$  & $u_1$ & $u_3$ \\
  $u_2$ & $0$ & $a$ & $u_2$ & $u_2$ & $c$ & $d$ & $1$ & $u_1$ & $u_2$ & $u_3$ \\
  $u_3$ & $u_3$ & $u_3$ & $u_3$ & $u_3$ & $u_3$ & $d$ & $1$ & $u_3$ & $u_3$ & $u_3$ \\ \hline
\end{tabular}
\caption{The internal uninorm $U_2$ in Example~\ref{Exa-1}.}
\label{tab-IntUni-2}
\end{table*}
\end{exa}

\begin{prot}\label{prot-ConsAuBuMore}
Suppose that $u\wedge e=v\wedge e$ and $u\vee e= v\vee e$ for all $u, v\in I_e$. Consider $u_0\in I_e$. Then the following results hold:
\begin{enumerate}[\rm(i)]
\item $\{\mathcal A_u\}_{u\in I_e}$ and $\{\mathcal B_u\}_{u\in I_e}$ fulfill {\rm Proposition~\ref{prot-existenceU0}(ii)} and {\rm Proposition~\ref{prot-existenceQ}(ii)}, where $\mathcal A_{u_0}\subseteq [0, u_0\wedge e]$ is an arbitrary lower set, $\mathcal A_u=[0, u\wedge e]$ for all $u\in I_e\setminus\{u_0\}$ and $\mathcal B_u=[u\vee e, 1]$ for all $u\in I_e$;

\item $\{\mathcal A_u\}_{u\in I_e}$ and $\{\mathcal B_u\}_{u\in I_e}$ fulfill {\rm Proposition~\ref{prot-existenceU0}(ii)} and {\rm Proposition~\ref{prot-existenceQ}(ii)}, where $\mathcal B_{u_0}\subseteq [u_0\vee e, 1]$ is an arbitrary upper set, $\mathcal A_u=[0, u\wedge e]$ for all $u\in I_e$ and $\mathcal B_u=[u\vee e, 1]$ for all $u\in I_e\setminus\{u_0\}$.
\end{enumerate}
\end{prot}
\begin{proof}
The proof is omitted as an easy exercise.
\end{proof}

\begin{exa}
Consider the bounded lattice $L=\{0, a, b, e, f, 1, u, v\}$ shown in Figure~\ref{figure-2}. For any lower set $\mathcal A_u\subseteq [0, u\wedge e]$, $\mathcal A_v=[0, v\wedge e]$, and $\mathcal B_u=\mathcal B_v=[u\vee e, 1]$, due to Proposition~\ref{prot-ConsAuBuMore}(i), there exists an internal uninorm on $L$ with neutral element $e$. Similarly, due to Proposition~\ref{prot-ConsAuBuMore}(ii), different internal uninorms on $L$ with neutral element $e$ can be constructed.
\end{exa}

\section{Concluding remarks}\label{sec5}

Under the assumption that all elements incomparable with the neutral element (collected in the set $I_e$) constitute an antichain, the structure of internal uninorms with neutral element $e$ is characterized in Theorem~\ref{thrm-characterizationthrm}. Based on the characterization theorem, the existence of an internal uninorm with neutral element $e$ is characterized from the perspective of the structure of the underlying lattice (see Theorem~\ref{thrm-chainAndantichain}). Note that the existence of an internal uninorm with neutral element $e$ is characterized when $I_e$ is a chain in~\cite{OZWD24}. For future work, it is natural to completely solve the existence problem by considering the remaining case that $I_e$ is neither a chain nor an antichain.

\end{document}